\documentclass[10pt]{article}
\usepackage{amsmath}
\usepackage{amssymb}
\usepackage{theorem}
\usepackage{appendix}

\newtheorem{lemma}{Lemma}[section]

\newtheorem{theorem}[lemma]{Theorem}
\newtheorem{corollary}[lemma]{Corollary}

\newcommand{\R}{{\bf R}}

\newcommand{\BB}{{\bf B}}
\newcommand{\HH}{{H}}

\newcommand{\pend}{$\Box$\\ }

\newcommand{\lam}{\lambda}

\newcommand{\Ome}{\Omega}
\newcommand{\vecy}[1]{\boldsymbol{#1}}
\newcommand{\modulus}[1]{\arrowvert#1\arrowvert}
\newcommand{\mint}{\int\limits_\Ome}
\newcommand{\hint}{\int\limits_{\R^{N}}}

\newcommand{\lap}{{\Delta}}

\newcommand{\bra}[1]{\left( #1 \right)}
\newcommand{\sbra}[1]{\small( #1 \small)}
\newcommand{\Quad}[1]{{Q \bra{#1}}}
\newcommand{\tQuad}[1]{{Q_{\lam\psi} \bra{#1}}}

\newcommand{\norm}[1]{\Vert #1 \Vert}
\newcommand{\cof}{a_{\vecy{\alpha},\vecy{\beta}}\bra{x}}
\newcommand{\sg}{e^{-Ht}}

\renewcommand{\Re}{{\rm Re}\;}

\newcommand{\Proof}{\underbar{Proof}{\hskip 0.1in}\par}

\newcommand{\beq}{\begin{equation}}

\newcommand{\g}{{\tilde{g}}}

\newcommand{\mlap}[1]{\norm{\bra{-\lap}^{#1}f}}

\newcommand{\kik}[1]{C^{\infty}_{c}\bra{#1}}
\newcommand{\ip}[2]{\langle #1,#2\rangle}

\newcommand{\twist}{e^{\lam \psi}}

\newcommand{\itwist}{e^{-\lam \psi}}

\newcommand{\sspace}{\,\,}

\newcommand{\colv}[2]{\begin{pmatrix} #1\\ #2 \end{pmatrix}}

\newcommand{\sobz}{W^{m,2}_0}
\newcommand{\dx}{d_x}
\newcommand{\dy}{d_y}
\newcommand{\ft}{\hat{f}\bra{\xi}}
\newcommand{\mft}{|\ft|}

\newcommand{\ggamma}{{m\bra{1-\epsilon}-\tfrac{N}{2}}}

\begin{document}
\title{\bf{Gaussian Upper Bounds on Heat Kernels of Uniformly Elliptic Operators on
Bounded Domains}}
\author{Narinder S Claire}
\date{}

\maketitle

\begin{abstract}
We obtain Gaussian upper bounds for heat kernels of higher order differential operators with Dirichlet boundary conditions on 
bounded domains in $\R^N$. The bounds  exhibit explicitly the nature of the spatial decay of the heat kernel close to the boundary as well as 
the long-time exponential decay implied by the spectral gap.  We make no smoothness 
assumptions on our operator coefficients which we assume only to be bounded and measurable.\\ \\
{\bf AMS Subject Classification : 35K25}  
\end{abstract} 

\section*{Introduction}

Off-diagonal Gaussian upper bounds for higher order differential operators with bounded measurable coefficients were first obtained
by Davies \cite{Da1}. The operators considered in that exposition were of order $2m$ on $L^2\bra{\R^N}$ with 
$2m >N$. The higher order operators, expressed as:
\begin{equation}
 Hf\bra{\vecy{x}}\quad :=\quad \sum\limits_{\substack{|\vecy{\alpha}|\leq m \\ |\vecy{\beta}|\leq m}}
 \bra{-1}^{|\vecy{\alpha}|}
D^{\vecy{\alpha}}\bra{a_{\vecy{\alpha},\vecy{\beta}}\bra{\vecy{x}}D^{\vecy{\beta}}f\bra{\vecy{x}}}
\end{equation}
were shown to have heat kernels with off-diagonal bounds demonstrated in the inequality:
\begin{equation}\label{gaussian}
|k\bra{t,x,y}| \leq
c_1 t^{-\frac{N}{2m}}
\exp \bra{-c_2 \frac{|x-y|^{\frac{2m}{2m-1} }}{t^\frac{1}{2m-1}}+c_3t}
\end{equation}
Subsequently Barbatis and Davies \cite{BD} were able to obtain optimal values for the constants $c_2$
and $c_3$ in terms of the ellipticity ratio and dimension.
\par
In this paper we address the question of upper bounds on heat kernels generated by uniformly elliptic differential operators 
with Dirichlet boundary conditions on bounded regions of $\R^N$. We make the same assumptions of the coefficients as in \cite{Da1}, namely
that they are measurable and bounded, and consequently find it more convenient to carry out the analysis with greater focus on the corresponding
quadratic forms. We do however assume that the quadratic forms, on a bounded region $\Ome\subset\R^N$, satisfy the ellipticity condition: 
\begin{equation}\label{H1}
c^{-1}\mlap{\frac{m}{2}}^2_2\leq\Quad{f}\leq c\mlap{\frac{m}{2}}^2_2
\end{equation}
for each $f \in \kik{\Ome}$ and $c$ strictly positive.
\par
Although  \eqref{gaussian} holds for operators on bounded regions with Dirichlet boundary conditions, the bounds do not reflect either the
spatial decay near the boundary or the long-time asymptotics. We extract the manner in which  
$$k\bra{t,x,y}
\rightarrow 0 \quad \text{as } x,y\rightarrow \partial\Ome $$ and show that the heat kernel has off-diagonal bounds demonstrated in the inequality:
$$|k\bra{t,x,y}| \leq c_1\sbra{1-\tfrac{N+2\gamma}{m}}^{-1}\,t^{-\frac{N+2\gamma}{2m}}\,d\sbra{x,\partial\Ome}^{\gamma}
\,d\sbra{y,\partial\Ome}^{\gamma}\exp\bra{-c_2 \frac{\modulus{x-y}^{\frac{2m}{2m-1}}}{t^{\frac{1}{2m-1}}}-st}$$
where $0\leq\gamma<m-\tfrac{N}{2}$ and $s$ is the spectral gap.\par 
The techniques we employ are close to those employed in \cite{Da1}, but we give more emphasis to the analysis
of spatial derivatives of the heat kernel. Moreover the exponential time decay is deduced by exploiting the spectral gap.
\par
Sharp off diagonal heat kernel bounds were also obtained by Barbatis \cite{Ba}
in terms of a non-euclidean metric based on the coefficients of
the operator, replacing the term $|x-y|$ in \eqref{gaussian} by $d\bra{x,y}$. In the case of bounded regions, an off
diagonal bound was obtained for highly non-convex regions for the
uniformly elliptic operator by Owen \cite{Ow}, in which he used the geodesic
distance but boundary behaviour was not the focus of that analysis. \par
Throughout this paper we will assume that $c$ and $c_i$ represent strictly positive constants. 
 
\subsection*{Notation}

Given $\vecy{\alpha}\in \R^N$ representing the multi-index $\bra{\alpha_1,\alpha_2,\alpha_3,...} $
where $$|\vecy{\alpha}|:=\alpha_1 +\alpha_2 +\alpha_3 +.....+\alpha_N$$
 we define the corresponding operator:
$$D^{\vecy{\alpha}}f :=\frac{\partial^{|\vecy{\alpha}|}f}{\partial^{\alpha_1}x_1
\partial^{\alpha_2}x_2
\partial^{\alpha_3}x_3....\partial^{\alpha_N}x_N}$$
and the set $V_{\vecy{\alpha}}$ such that  $$V_{\vecy{\alpha}} =\{\vecy{r} : r_i\leq\alpha_i
\text{ for all } i \}$$
Moreover given any multi-index $\vecy{r}$ in $V_{\vecy{\alpha}}$
we define the vector factorial as:
$$\colv{\vecy{\alpha}}{\vecy{r}}:= \colv{\alpha_1}{r_1}\colv{\alpha_2}{r_2}\colv{\alpha_3}{r_3}...\colv{\alpha_N}{r_N}$$ \\
The directional derivative of order $m$ of an appropriately smooth function along a vector $\vecy{v}$ in $\R^N$ is expressed as
$$ \frac{\partial^m f}{\partial \vecy{v}^m}:=\cdots\nabla\bra{\nabla\bra{\nabla f . \vecy{v}}.\vecy{v}}.
\vecy{v}$$\par

\subsection*{Quadratic Form}
It is helpful to give an indicative though a non-rigorous formulation of the family of higher order operators that we focus on 
in this the paper. The operator is defined more completely through it's quadratic form.
Given a bounded domain $\Ome$ in $\R^N$ we express the operator of order $2m>N$ as : 
\begin{equation}
 Hf\bra{\vecy{x}}\quad :=\quad \sum\limits_{\substack{|\vecy{\alpha}|\leq m \\ |\vecy{\beta}|\leq m}}
 \bra{-1}^{|\vecy{\alpha}|}
D^{\vecy{\alpha}}\bra{a_{\vecy{\alpha},\vecy{\beta}}\bra{\vecy{x}}D^{\vecy{\beta}}f\bra{\vecy{x}}}
\end{equation}
where $a_{\vecy{\alpha},\vecy{\beta}}$ are complex bounded measurable functions.\\ 
The associated quadratic form $Q$ 
\begin{equation}
Q\bra{f}\quad :=\quad \sum\limits_{\substack{|\vecy{\alpha}|\leq m \\ |\vecy{\beta}|\leq m}}
\int_{\substack{\Ome}}a_{\vecy{\alpha},\vecy{\beta}}\bra{\vecy{x}}D^{\vecy{\beta}}f\bra{\vecy{x}}\overline{D^{\vecy{\alpha}}f\bra{\vecy{x}}}
\end{equation}

defined with domain equal to the Sobolev space 
$\sobz \bra{\Ome}$ will be assumed to satisfy the ellipticity condition with a strictly positive constant $c$ 
\begin{equation}\label{H1}
c^{-1}\mlap{\frac{m}{2}}^2_2\leq\Quad{f}\leq c\mlap{\frac{m}{2}}^2_2
\end{equation}
We define the the spectral gap
$$ s=\inf\limits_{f\in \kik{\Ome}}\frac{\Quad{f}}{\norm{f}_2^2}$$
Since we have made the assumption that $\tfrac{N}{2m}<1$, it will be informative to track the dependency on this constraint by defining the 
quantity $0<\epsilon\leq1-\tfrac{N}{2m}$ and 
$$\gamma:=\ggamma$$ \par
For a given point $x$ in $\Ome$ we define it's distance from the boundary $\partial \Ome$ 
$$\dx:= d\bra{\vecy{x},\partial\Ome} = \inf\limits_{\vecy{y}\in \partial\Ome}\modulus{\vecy{x} - \vecy{y}}$$
We define the function $\g\bra{t}$ such that 
\begin{equation}\label{e:gt}
\g\bra{t} := \begin{cases} s e^{-2s t} &  t>\frac{1}{s}\\
\frac{1}{t}e^{-s t-1} &  t \leq \frac{1}{s}\\
\end{cases}
\end{equation}


\section{Boundary Behaviour}
Having imposed Dirichlet boundary conditions on our operator we expect the heat kernel $k\bra{t,x,y}$ to vanish at the boundary. The
precise nature of this decay can be deduced by application of the Sobolev embedding theorem : $$\sobz \bra{\Ome} \hookrightarrow
C^{\gamma}_0\bra{\Ome}$$
We consider the norm 
on $C^{\gamma}_0\bra{\Ome}$ to be defined as 
$$\norm{f}_{C^{\gamma}_0} := \sum\limits_{|\alpha| \leq n} \norm{D^{\alpha}f}_{\infty} \:
+ \sup\limits_{\vecy{x},\vecy{y} \in \Ome}\frac{|D^{\alpha}f\bra{\vecy{y}}-D^{\alpha}f\bra{\vecy{x}}|}
{|\bra{\vecy{y}-\vecy{x}}|^{\kappa}}$$
where $n$ and $\kappa$ are the integer and the fractional parts of $\gamma$ respectively.

\begin{lemma}\label{sinbad}
There is a strictly positive constant $c$ such that for all $f \in \sobz \bra{\Ome}$ and any unit vector $\vecy{v}$ in $\R^N$
\begin{equation}\label{e:premature}
\modulus{\tfrac{\partial^n f\bra{\vecy{x}}}{\partial \vecy{v}^n}}\leq \tfrac{c}{\sqrt{\epsilon}}\sspace
{\dx^{\kappa}}\sspace\Quad{f}^{\frac{1-\epsilon}{2}}\sspace\norm{f}_2^{\epsilon}
\end{equation}
for all $x \in \Ome$ 
\end{lemma}
\Proof
By applying Fourier transform to $\frac{\partial^nf}{\partial \vecy{v}^n}$
we have for all $x$ in $\Ome$
$$\sup\limits_{\vecy{y} \in \R^N}\frac{\modulus{\tfrac{\partial^n f\bra{\vecy{x}}}{\partial \vecy{v}^n}
-\frac{\partial^n f\bra{\vecy{y}}}{\partial \vecy{v}^n}   }}
{|\bra{\vecy{y}-\vecy{x}}|^{\kappa}} \leq c\hint |\vecy{\xi}|^{\vecy{\gamma}}\, |\ft|d \vecy{\xi}$$
and consequently 
\begin{eqnarray*}
\modulus{\tfrac{\partial^n f\bra{x}}{\partial \vecy{v}^n} } \leq c\, \dx^{\kappa} \hint |\vecy{\xi}|^{\vecy{\gamma}}\, |\ft|d \vecy{\xi}
\end{eqnarray*}

hence for a positive $\mu$
$$
\modulus{\tfrac{\partial^n f\bra{x}}{\partial \vecy{v}^n}}^2
\leq c \dx^{2\kappa} \bra{\hint\bra{\mu + |\vecy{\xi}|^2}^{\frac{\gamma-m}{2}}\bra{\mu + |\vecy{\xi}|^2}^{\frac{m}{2}}\mft d^N \vecy{\xi} }^2   \\
$$
It then follows from Cauchy-Schwartz that 

\begin{eqnarray*}
\modulus{\tfrac{\partial^n f\bra{x}}{\partial \vecy{v}^n}}^2&\leq & c\, \dx^{2\kappa}\, \bra{\hint\bra{\mu + |\vecy{\xi}|^2}^m \mft^2 d^N \vecy{\xi} }
\bra{\hint \frac{1}{\bra{\mu + |\vecy{\xi}|^2}^{m-\gamma}} d^N \vecy{\xi} }\\
& \leq & \frac{c}{\epsilon}\, \dx^{2\kappa}\,\bra{\mu^{m\bra{1-\epsilon}}\norm{f}^2_2 +\mu^{-m\epsilon}\Quad{f}}
\end{eqnarray*}
Optimizing over $\mu$ to find that 
$$\mu^m = \frac{\epsilon}{1-\epsilon}
\frac{\Quad{f}}{\norm{f}_2^2}$$
completes the proof.\\
\pend
We proceed to find an upper bound for the heat kernel by applying lemma \ref{sinbad} to $f_t:=\sg f$ but first we need a more comprehensive upper bound
for  $\modulus{\tfrac{\partial^n f_t\bra{x}}{\partial \vecy{v}^n}}$ by applying the Spectral Theorem.

\begin{lemma}\label{SGsinbad}
If $f_t$  is $ \sg f$ for some $f$ in $L^2\bra{\Ome}$ then 
\begin{equation}
\Quad{f_t} \leq \g\bra{t} \,\norm{f}^2_2
\end{equation}

\end{lemma}
\Proof
The inequality follows from 
$$\Quad{f_t} \leq  \norm{H \sg}\,\norm{\sg}\,\norm{f}^2_2$$
and an application of the Spectral Theorem.\\
\pend

We can now combine lemmas \ref{sinbad} and  \ref{SGsinbad} to yield our upper bound for the heat kernel.

\begin{lemma}\label{lem:sun}
The heat kernel $k\bra{t,x,y}$ generated by the differential operator $H$ satisfies the inequalities
\begin{eqnarray*}
\modulus{k\sbra{t,x,y}}\leq c\sspace{\small(1-\tfrac{N+2\gamma}{2m}\small)}^{-1}\sspace t^{-\frac{N+2\gamma}{2m}}\sspace \dx^{\gamma}\sspace\dy^{\gamma} & {\text{when }} t<\frac{2}{s}\\
\\
\modulus{k\sbra{t,x,y}}\leq c\sspace{\small(1-\tfrac{N+2\gamma}{2m}\small)}^{-1}\sspace e^{-st}\sspace\dx^{\gamma}\sspace\dy^{\gamma}  &  {\text{when }} t\geq\frac{2}{s}
\end{eqnarray*}

\end{lemma}
\Proof
If $f_t:=\sg f$ then from lemma \ref{sinbad} we have  
\begin{eqnarray*}\label{e:premature}
\modulus{\tfrac{\partial^n f_t\bra{\vecy{x}}}{\partial \vecy{v}^n}}\leq \tfrac{c}{\sqrt{\epsilon}}\sspace
{\dx^{\kappa}}\sspace\Quad{f_t}^{\frac{1-\epsilon}{2}}\sspace\norm{f_t}_2^{\epsilon}
\end{eqnarray*}
Choosing $\vecy{v}$ appropriately and integrating yields
\begin{eqnarray*}
\modulus{\mint k\bra{t,x,u}f\bra{u} \ d^N u}
&\leq 
\tfrac{c}{\sqrt{\epsilon}}\sspace {\dx^{\gamma}}\sspace\g\bra{t}^{\frac{1-\epsilon}{2}}
\sspace\norm{f}_2^{1-\epsilon}\sspace\norm{f_t}_2^{\epsilon}\\
&\leq
\tfrac{c}{\sqrt{\epsilon}}\sspace{\dx^{\gamma}}\sspace\g\bra{t}^{\frac{1-\epsilon}{2}}e^{-st\epsilon }\sspace\norm{f}_2 
\end{eqnarray*}
This inequality gives us a bound for the $L_2$ norm of $k\bra{t,y,\cdot}$
$$ \norm{k\bra{t,y,\cdot}}_2 \leq \tfrac{c}{\sqrt{\epsilon}}\sspace{\dx^{\gamma}}\sspace\g\bra{t}^{\frac{1-\epsilon}{2}}e^{-st\epsilon }$$
which can then be applied to obtain the required bounds
\begin{eqnarray*}
|k\bra{t,x,y}|&=&|\mint k\bra{\tfrac{t}{2},x,u}k\bra{\tfrac{t}{2},u,y}\ d^N u |\\&\leq&
 \norm{k\bra{\tfrac{t}{2},x,\cdot} }_2 \norm{k\bra{\tfrac{t}{2},y,\cdot} }_2  \\
&\leq &
\tfrac{c}{\epsilon}\sspace{\dx^{\gamma}}\sspace{\dy^{\gamma}}\sspace\g\bra{\tfrac{t}{2}}^{1-\epsilon}{e}^{-st\epsilon }
\end{eqnarray*}
\pend

\section{Gaussian Bounds}

Gaussian bounds exhibiting boundary decay can be given by interpolation between the bounds found in lemma \ref{lem:sun}
and those obtained in Davies \cite{Da1}.The drawback of this method however, is that the presence of the term $c_3t$ does not
imply the the long-time exponential decay that we expect and would like to show.

Dirichlet boundary conditions imply a positive spectral gap and hence exponential time {\em decay}. We give bounds that 
more concisely exhibit this behaviour.
\par One of the key features of hypothesis \eqref{H1} is that many required operator inequalities can be reduced to  proving the corresponding inequalities
for polynomial symbols by applying Fourier transforms.

Barbatis and Davies \cite{BD} make a stronger assumption on the operator coefficients.
They assume that for a strictly positive $q$

\begin{equation}
 q^{-1}\sum\limits_{\substack{|\vecy{\alpha}|= m \\ |\vecy{\beta}|= m}}
 a_{0,\vecy{\alpha}\vecy{\beta}}\xi_{\vecy{\alpha}}\overline{\xi_{\vecy{\beta}}} \leq
 \sum\limits_{\substack{|\vecy{\alpha}|= m \\ |\vecy{\beta}|= m}}
 a_{\vecy{\alpha}\vecy{\beta}}\bra{x}\xi_{\vecy{\alpha}}\overline{\xi_{\vecy{\beta}}}\leq
 q \sum\limits_{\substack{|\vecy{\alpha}|= m \\ |\vecy{\beta}|= m}}
 a_{0,\vecy{\alpha}\vecy{\beta}}\xi_{\vecy{\alpha}}\overline{\xi_{\vecy{\beta}}}
 \end{equation}
 almost everywhere in $\Ome$ and the non-negative
 coefficient matrix $\bf{A}_0 =  a_{0,\vecy{\alpha}\vecy{\beta}}$ satisfies 
 \begin{equation}
 \ip{\bra{-\lap}^m f}{g}=\mint  a_{0,\vecy{\alpha}\vecy{\beta}}D^{\vecy{\alpha}}f\overline{D^{\vecy{\beta}}g} d^N \vecy{x}
 \end{equation}
 for all functions $f,g \in \kik{\Ome}$.

They obtained the necessary estimates for the polyharmonic operator by way of Fourier
transforms and an application of the polarization identity on the co-efficient
matrix. We do not make this assumption here.

\subsection{Twisted Quadratic Form Inequality}

Given  $\vecy{x_0}$ in $\Ome$ and a unit vector $\vecy{a}$ in $\R^N$, we 
define the bounded function $\psi_{\vecy{x_0},\vecy{a}}$ on $\Ome$
\begin{equation}\label{psipsi}
\psi_{\vecy{x_0},\vecy{a}}\bra{\vecy{x}}:=\langle \vecy{x}-\vecy{x_0},\vecy{a}\rangle
\end{equation}
Moreover given a real number $\lambda$ and dropping the subscripts on $\psi_{\vecy{x_0},\vecy{a}}$ 
we define the multiplicative operator $e^{\lambda\psi}$
$$ e^{\lambda\psi}:L^2\bra{\Ome}\rightarrow L^2\bra{\Ome}$$
$$ \bra{e^{\lambda\psi}f}\bra{x} := e^{\lambda\psi\bra{x}}f\bra{x}$$
$e^{\lambda\psi}$ is a homeomorphism on $L^2\bra{\Ome}$.
\par
The twisted operator $ \HH_{\lambda\psi}$ is then defined as follows:
$$ \HH_{\lambda\psi}f:=\itwist \HH \twist f$$
with
$$ Dom\bra{\HH_{\lambda\psi}}=\{\itwist f : f\in Dom\bra{\HH}\}$$
More importantly we define the twisted quadratic form, $Q_{\lam \psi}$ as follows:
\begin{eqnarray*}
Q_{\lam \psi}\bra{f}& :=&\Quad{\twist f, \itwist f}\\
&=&\sum_{\vecy{\alpha} ,\vecy{\beta}}\quad \mint
\cof \bra{D^{\vecy{\alpha}}\bra{\twist f}}\bra{D^{\vecy{\beta}}\bra{\itwist \overline{f}}}
 d^{N}\boldsymbol{x}\\
&=&\sum_{\vecy{\alpha} ,\vecy{\beta}}\quad \mint
\cof \bra{\itwist D^{\vecy{\alpha}}\bra{\twist f}}\bra{\twist D^{\vecy{\beta}}\bra{\itwist \overline{f}}}
d^{N}\boldsymbol{x}\\
\end{eqnarray*}

By Leibniz, for each  $ f \in \kik{\Ome}$ we can expand
\begin{equation}
\itwist D^{\vecy{\alpha}}\twist f = \sum_{\vecy{r} \in V_{\vecy{\alpha}}}
\colv{\vecy{\alpha}}{\vecy{r}} \lam^{|\vecy{\alpha}-\vecy{r}|}
\vecy{a}^{\vecy{\alpha}-\vecy{r}}D^{\vecy{r}} f
\end{equation}
\par
It is then possible to show that the difference $per\bra{\lam} := Q_{\lam \psi}\bra{f}-Q \bra{f}$ is 
$$
\sum_{\vecy{\alpha} ,\vecy{\beta}} \mint
\cof \sum\nolimits '\colv{\vecy{\alpha}}{\vecy{r}}\colv{\vecy{\beta}}{\vecy{s}}
\lam^{|\bra{\vecy{\alpha}+\vecy{\beta}}-\bra{\vecy{r}+\vecy{s}}|}
\vecy{a}^{\vecy{\alpha}-\vecy{r}}\vecy{a}^{\vecy{\beta}-\vecy{s}}
D^{\vecy{r}} f D^{\vecy{s}} \overline{f}
$$
where the summation $\sum\nolimits '$ runs over all $\vecy{r}$ in $V_{\vecy{\alpha}}$ and $\vecy{s}$ in $V_{\vecy{\beta}}$ but where either
 $|\vecy{r}|< |\vecy{\alpha}|$ or $|\vecy{s}|<|\vecy{\beta}|$. The terms where both $\vecy{r}=\vecy{\alpha}$
and $\vecy{s}=\vecy{\beta}$  are incorporated in  $\Quad{f}$

Since the coefficients of the operator are uniformly bounded we have for $c>0$
$$
|per(\lam)|\leq c \sum_{\vecy{\alpha},\vecy{\beta}}\sum\nolimits '|\lam|^{|\bra{\vecy{\alpha}+\vecy{\beta}}-\bra{\vecy{r}+\vecy{s}}|}
\mint |D^{\vecy{r}}f||D^{\vecy{s}}f| d^{N} \vecy{x} \\
$$
and applying Cauchy-Schwartz 
$$
|per(\lam)|\leq  \sum_{1 \leq |\vecy{\alpha}|=|\vecy{\beta}|\leq m}\sum\nolimits '
\norm{\lam^{|\vecy{\alpha}|-|\vecy{r}|}D^{\vecy{r}}f}_{2}
\norm{\lam^{|\vecy{\beta} |-|\vecy{s}|}D^{\vecy{s}}f}_{2}
$$

From lemma \ref{lem:epsilon} we see that for $\varepsilon<1$
\begin{eqnarray*}
|per(\lam)| &\leq & c \sum_{1\leq p \leq m} \sum \nolimits '\norm{|\lam|^{p-|\vecy{r}|}\bra{-\lap}^{\frac{|\vecy{r}|}{2}}f}_{2}
\norm{\lam^{p-|\vecy{s}|}\bra{-\lap}^{\frac{|\vecy{s}|}{2}}f}_{2}\\
&\leq & c\sum_{1\leq p \leq m}
\varepsilon \mlap{\frac{p}{2}}^{2}_{2}+\varepsilon^{1-2p}\lam^{2p}\norm{f}^{2}_{2}
\end{eqnarray*}
It follows from lemma \ref{lem:stephen} with $\rho=\varepsilon^{1-2m}$
$$|per(\lam)|\leq  c_1\varepsilon\bra{1+\theta}\Quad{f} + c_1\varepsilon^{1-2m} \bra{1+\theta s\epsilon^{2m-1}}^{2m} \lambda^{2m}\norm{f}^{2}_{2}$$
and simplifying
$$|per(\lam)|\leq  c_1\varepsilon\bra{1+\theta}\Quad{f} + c_1\varepsilon^{1-2m} \bra{\left[1+\theta s\right] \lambda}^{2m}\norm{f}^{2}_{2}$$

Using 
$$Re\left[Q \bra{f}-Q_{\lam \psi}\bra{f}\right]  = Q \bra{f}-Re\left[Q_{\lam \psi}\bra{f}\right] \leq \modulus{Q \bra{f}-Q_{\lam \psi}\bra{f}}$$
we have
\begin{equation}\label{rea}
 \left[1-c_1\varepsilon\bra{1+\theta}\right]\Quad{f} \leq  Re {\tQuad{f}} + c_1\varepsilon^{1-2m} \bra{\left[1+\theta s\right] \lambda}^{2m}\norm{f}^{2}_{2}
\end{equation} 
Moreover 
\begin{equation}\label{img}
|Im {\tQuad{f}}| \leq c_1\varepsilon\bra{1+\theta}\Quad{f} + c_1\varepsilon^{1-2m} \bra{\left[1+\theta s\right] \lambda}^{2m}\norm{f}^{2}_{2}
\end{equation}
We note from lemma \ref{lem:stephen} that these inequalities are true for all positive $\theta$ 
\subsection{Spectral Gap}

From the outset we know that $\HH$ has a positive definite least eigenvalue $s$. We define 
$$ \hat{H} : = \HH -s $$
similarly
\begin{equation}\label{zerohat} \hat{\HH}_{\lambda\psi}f:=\itwist \HH \twist f\end{equation} 
and crucially
$$\hat{Q}_{\lam \psi}\bra{f} := Q_{\lam \psi}\bra{f} -s\norm{f}_2^2$$
It is easy to see that the real an imaginary parts of the newly defined twisted form, $\hat{Q}_{\lam \psi}$ satisfy
\begin{equation}\label{reaprime}
 \left[1-c_1\epsilon\bra{1+\theta}\right]\hat{Q}\bra{f} \leq  Re {\hat{Q}_{\lam \psi}\bra{f}} + c_1\epsilon^{1-2m} \bra{\left[1+\theta s\right] \lambda}^{2m}\norm{f}^{2}_{2}
\end{equation} 
and
\begin{equation}\label{imgprime}
|Im {\hat{Q}_{\lam \psi}\bra{f}}| \leq c_1\epsilon\bra{1+\theta}\hat{Q}\bra{f} + c_1\epsilon^{1-2m} \bra{\left[1+\theta s\right] \lambda}^{2m}\norm{f}^{2}_{2}
\end{equation}

\subsection{Twisted Semigroup Inequalities}
\begin{lemma}\label{p}
There is a positive constant $k$ such that for all $c\geq k$ and all $0<p<1$

$$Q_{\lambda,\psi}'\bra{f} := \hat{Q}_{\lam \psi}\bra{f} + c\bra{1+p}\bra{1+s}^{2m}\lambda^{2m}\norm{f}^{2}_{2}$$
is a sectorial form with vertex at $0$ and a semi-angle smaller then $\frac{\pi}{2}$.\\
Moreover 
$$
 \frac{p}{1+p}\hat{Q}\bra{f}\leq Re {\hat{Q}_{\lam \psi}\bra{f}} + c\bra{1+p}\bra{1+s}^{2m}\lambda^{2m}\norm{f}^2_2
$$
\end{lemma}
\Proof
With $\theta=\frac{2-\bra{1+p}}{\bra{1+p}}$ and $\varepsilon=\frac{1}{2c_1}$ the second part of the RHS of (\ref{imgprime})
$$
c_1\varepsilon^{1-2m} \bra{\left[1+\theta s\right] \lambda}^{2m}\norm{f}^{2}_{2}
$$
is less than 
$$
c_1^{2m}2^{2m-1}\bra{1+s}^{2m}\lambda^{2m}\norm{f}^{2}_{2}
$$
hence
\begin{eqnarray*}
p|Im Q_{\lambda,\psi}'\bra{f}| &\leq& \frac{p}{1+p}\hat{Q}\bra{f}+ pc_1^{2m}2^{2m-1}\bra{1+s}^{2m}\lambda^{2m}\norm{f}^{2}_{2}\\
&=& \left[1-c_1\epsilon\bra{1+\theta}\right]\hat{Q}{f} + pc_1^{2m}2^{2m-1}\bra{1+s}^{2m}\lambda^{2m}\norm{f}^{2}_{2}\\
&\leq& Re {\hat{Q}_{\lam \psi}\bra{f}} +\bra{1+p}\bra{c_1^{2m}2^{2m-1}\bra{1+s}^{2m}\lambda^{2m}}\norm{f}^{2}_{2}\\
&=& Re {\hat{Q}_{\lam \psi}\bra{f}} + c_2\bra{1+p}\bra{1+s}^{2m}\lambda^{2m}\norm{f}^2_2\\
&=&Re {Q_{\lambda,\psi}'\bra{f}}
\end{eqnarray*}
This then implies that
$$\frac{|Im Q_{\lambda,\psi}'{f}|}{Re Q_{\lambda,\psi}'{f}}\leq\frac{1}{p}$$
and consequently $$ |Arg \, Q_{\lambda,\psi}'\bra{f}|\leq\tan^{-1}\bra{\frac{1}{p}}<\frac{\pi}{2}$$
\pend
We define the operator $ H_{\lambda,\psi}'$ 
$$ H_{\lambda,\psi}' :=  \hat{H}_{\lambda,\psi}+c\bra{1+s}^{2m}\lambda^{2m}$$
where $c>0$ such that 
$$ |Arg\ip{H_{\lambda,\psi}'f}{f}|<\frac{\pi}{2}$$
We recall the following corollary from Kato\cite{ka}.
\begin{corollary} \label{his}
Let $\beta$ be angle such that $\frac{\pi}{4}<\beta<\frac{\pi}{2}$ and set 
$$p=\frac{1}{\tan \beta}$$
Then $e^{-H_{\lambda,\psi}'z}$ is an analytic semigroup in the sector
$$S_{\beta}=\{z:|Arg\bra{z}|<\beta\}$$
moreover $$\norm{e^{-H_{\lambda,\psi}'z}}\leq 1$$
\end{corollary}
\Proof
see[p492:\cite{ka}]\\ \pend
From  Davies \cite{Da6} Theorem 2.38 it is evident that there is a positive constant $c_p$ such that
\begin{equation}\label{daviesanalytic}
\norm{H_{\lambda,\psi}'e^{-H_{\lambda,\psi}'z}}\leq \frac{c_p}{Re\,z}
\end{equation}
for all $z\in S_\beta$.\\
In the following estimates we let $\hat{H}_{\lam \psi}$ be defined by (\ref{zerohat}).
\begin{lemma}\label{hisanalytic}
There is a positive constant $c$ such that 
\begin{equation}\label{swappy}
\norm{e^{-\hat{H}_{\lam \psi}t}}\leq \exp \big[c\bra{1+s}^{2m}\lambda^{2m}t\big]
\end{equation}
\end{lemma}
\Proof
Let $f \in L^2 \bra{\Ome}$ and define $f_t  :=e^{-\hat{H}_{\lam \psi}t} f$, then solving
\begin{eqnarray*}
\frac{d}{dt}\norm{f_t}_2^2 & = & \ip{-\hat{H}_{\lam \psi}f_t}{f_t}+
\ip{f_t}{-\hat{H}_{\lam \psi}f_t}\\
& = & -2Re \hat{Q}_{\lam \psi}\bra{f_t}\\
&\leq &  2c\bra{1+s}^{2m}\lambda^{2m}\norm{f_t}^2_2
\end{eqnarray*}
proves the claim.\\
\pend 
\begin{lemma} \label{interpolate}
Whenever $\beta\geq 0$, there is a positive constant $c_2$ such for any $0<\alpha<1$ we have 
\begin{equation}\label{mixed}
\norm{\hat{H}_{\lam \psi}e^{-\hat{H}_{\lam \psi}t}} + \beta \bra{1+s}^{2m}\lambda^{2m}\,\norm{e^{-\hat{H}_{\lam \psi}t}}\leq 
\frac{c_2}{\alpha t}e^{c\bra{1+\alpha}\bra{1+s}^{2m}\lambda^{2m}t}
\end{equation}
\end{lemma}
\Proof
Applying the triangle inequality gives us
$$ \norm{\hat{H}_{\lam \psi}e^{-\hat{H}_{\lam \psi}t}}e^{-c\bra{1+s}^{2m}\lambda^{2m}t} \leq 
\norm{H'_{\lambda,\psi}e^{-H'_{\lam \psi}t}} + c\bra{1+s}^{2m}\lambda^{2m}\norm{e^{-H'_{\lam \psi}t}}
$$
and it follows from (\ref{daviesanalytic}) and lemma \ref{hisanalytic} 
$$ \norm{\hat{H}_{\lam \psi}e^{-\hat{H}_{\lam \psi}t}} \leq 
\frac{c_1}{t}e^{c\bra{1+s}^{2m}\lambda^{2m}t}+ c\bra{1+s}^{2m}\lambda^{2m}e^{c\bra{1+s}^{2m}\lambda^{2m}t}
$$
Then the LHS of (\ref{mixed}) is bounded above by 
\begin{equation}\label{ohno}
\frac{c_1}{t}e^{c\bra{1+s}^{2m}\lambda^{2m}t}+ c\bra{1+\beta}\bra{1+s}^{2m}\lambda^{2m}e^{c\bra{1+s}^{2m}\lambda^{2m}t}
\end{equation}
For positive real $x$ we can re-write $xe^{xt}$ as 
$$\frac{1}{\alpha} e^{\bra{1+\alpha}xt}. x\alpha e^{-\alpha xt}$$
 for some $\alpha>0$.
We can then observe that $$xe^{xt} \leq \frac{e^{\bra{1+\alpha}xt}}{\alpha} . \frac{e^{-1}}{t}  $$ 
and see that similarly 
$$ 
c\bra{1+s}^{2m}\lambda^{2m}e^{c\bra{1+s}^{2m}\lambda^{2m}t} \leq \frac{e^{-1}}{\alpha t}e^{\bra{1+\alpha}c\bra{1+s}^{2m}\lambda^{2m}t}
$$
Substitution into (\ref{ohno}) attains the claimed inequality. \\
\pend
\begin{lemma}
Let $f \in L^2 \bra{\Ome}$ and define $\hat{f}_t \, :=e^{-\hat{H}_{\lam \psi}t} f$ then there are positive constants 
$c_1$ and $c_2$ such that
\begin{equation}\label{e:quad}
\hat{Q}\sbra{\hat{f}_t} \leq \frac{c_1}{\alpha t}e^{c_2\bra{1+s}^{2m}\lambda^{2m}t}\norm{f}_2^2
\end{equation}
\end{lemma}
\Proof
From lemma \ref{p} we have that for some positive constants $c_1$ and $c_2$ 
\begin{eqnarray*}
\hat{Q}\sbra{\hat{f}_t} &\leq& c_1\Re \hat{Q}_{\lam \psi}\sbra{\hat{f}_t}+c_2\bra{1+s}^{2m}\lambda^{2m}\norm{\hat{f}_t}_2^2\\
&\leq & c_1\norm{\hat{H}_{\lam \psi}\hat{f}_t}_2\,\norm{\hat{f}_t}_2+c_2\bra{1+s}^{2m}\lambda^{2m}\norm{\hat{f}_t}_2^2
\end{eqnarray*}
By using lemma \ref{interpolate} we can see that for $0<\alpha<1$
\begin{eqnarray*}
\hat{Q}\sbra{\hat{f}_t} &\leq& \frac{c_1}{\alpha t}e^{c\bra{1+\alpha}\bra{1+s}^{2m}\lambda^{2m}t}\norm{f}_2\,\norm{\hat{f}_t}_2\\
&\leq & \frac{c_1}{\alpha t}e^{2c\bra{1+\alpha}\bra{1+s}^{2m}\lambda^{2m}t}\norm{f}^2_2\\
\end{eqnarray*}
\pend
\begin{corollary}\label{quadinequality}
There are positive constants $c_1$ and $c_2$ such that 
 $$Q\sbra{e^{-H_{\lam \psi}t} {f}} \leq \frac{c_1}{\alpha t}e^{c_2\bra{1+s}^{2m}\lambda^{2m}t-2st}\norm{f}_2^2$$
\end{corollary}
\Proof
By Substitution from 
$$\hat{Q}\sbra{\hat{f}_t} = {Q}\sbra{\hat{f}_t} - s \norm{\hat{f}_t}_2^2$$
and 
$$\hat{f}_t =e^{-\hat{H}_{\lam \psi}t} f = e^{-H_{\lam \psi}t}e^{-st} f$$
it can be seen that 
$$\hat{Q}\sbra{\hat{f}_t} = e^{2st}Q\sbra{e^{-H_{\lam \psi}t} {f}} - s \norm{\hat{f}_t}_2^2$$
hence for $0<\alpha<1$
\begin{eqnarray*}
Q\sbra{e^{-H_{\lam \psi}t} {f}} &\leq& e^{-2st}\left[\hat{Q}\sbra{\hat{f}_t} + s\norm{\hat{f}_t}_2^2\right]\\
&\leq & e^{-2st}\left[\tfrac{c_1}{\alpha t}e^{c_2\bra{1+\alpha}\bra{1+s}^{2m}\lambda^{2m}t}\norm{f}_2^2 +
 s\norm{\hat{f}_t}_2^2\right] 
\end{eqnarray*}
Applying the estimate from lemma \ref{hisanalytic} completes the proof.\\
\pend


\subsection{Heat Kernel bounds}

\begin{theorem}
The integral kernel $k_{\lam,\psi}\bra{t,x,y}$ of $e^{-H_{\lam \psi}t}$ satisfies the inequality
\begin{equation}\label{twistedbound}
|k_{\lam,\psi}\bra{t,x,y}| \leq \frac{c}{\epsilon t^{1-\epsilon}}\,\dx^{\gamma}
\dy^{\gamma}e^{\small[c_2\bra{1+s}^{2m}\lambda^{2m}-s\small]t}\end{equation}
for some positive constants $c$ and $c_2$.
\end{theorem}
\Proof
For $f\in L^2\bra{\Ome}$ we define $f_t:=e^{-H_{\lam\psi}t}f$. From lemma \ref{sinbad}
$$\modulus{\tfrac{\partial^n f_t\bra{\vecy{x}}}{\partial \vecy{v}^n}}\leq \tfrac{c}{\sqrt{\epsilon}}\sspace
{\dx^{\kappa}}\sspace\Quad{f_t}^{\frac{1-\epsilon}{2}}\sspace\norm{f_t}_2^{\epsilon}
$$
and from corollary \ref{quadinequality}
\begin{equation}\label{aalmost}
\modulus{\tfrac{\partial^n f_t\bra{\vecy{x}}}{\partial \vecy{v}^n}}
\leq 
\tfrac{c}{t^{\frac{1-\epsilon}{2}}\sqrt{\epsilon} }\dx^{\kappa}
e^{\sbra{\frac{1-\epsilon}{2}}\small[c_2\bra{1+s}^{2m}\lambda^{2m}-2s\small]t}\norm{f}_2^{1-\epsilon}\,\norm{f_t}_2^{\epsilon}
\end{equation}

Recalling that $f_t=e^{-H_{\lam\psi}t}f=e^{-st}e^{-\hat{H}_{\lam\psi}t}f$ then with lemma \ref{hisanalytic}
we have the estimate
$$\norm{f_t}_2^{\epsilon}\leq  \exp \sbra{\epsilon c_2\sbra{1+s}^{2m}\lambda^{2m}t}e^{-\epsilon st}\norm{f}_2^\epsilon$$ 
Substituting this estimate into (\ref{aalmost})
\begin{equation*}
\modulus{\tfrac{\partial^n f_t\bra{\vecy{x}}}{\partial \vecy{v}^n}}
\leq 
\tfrac{c}{t^{\frac{1-\epsilon}{2}}\sqrt{\epsilon} }{\dx^{\kappa}}
e^{\small[c_2\bra{1+s}^{2m}\lambda^{2m}-s\small]t}\norm{f}_2
\end{equation*}

Integrating along the path to the boundary 
\begin{eqnarray*}
|\mint k_{\lam,\psi}\bra{t,x,u}f\bra{u} \ d^N u| 
 \leq
\tfrac{c}{t^{\frac{1-\epsilon}{2}}\sqrt{\epsilon} }{\dx^{\gamma}}
e^{\small[c_2\bra{1+s}^{2m}\lambda^{2m}-s\small]t}\norm{f}_2 
\end{eqnarray*}
Following a similar argument to that in the proof to lemma \ref{lem:sun}
\begin{eqnarray*}
|k_{\lam,\psi}\bra{t,x,y}|\leq \norm{k_{\lam,\psi}\bra{  \tfrac{t}{2}  ,x,\cdot}}_2
\norm{k_{\lam,\psi}\bra{  \tfrac{t}{2}   ,y,\cdot}}_2
\end{eqnarray*}
to yield the upper bound 
$$|k_{\lam,\psi}\bra{t,x,y}| \leq \frac{c}{\epsilon t^{1-\epsilon}}\dx^{\gamma}\dy^{\gamma}e^{\small[c_2\bra{1+s}^{2m}\lambda^{2m}-s\small]t}$$
\pend

\begin{theorem}
The integral kernel $k\bra{t,x,y}$ of $e^{-Ht}$ satisfies the inequality
$$|k\bra{t,x,y}| \leq c_1\sbra{1-\tfrac{N+2\gamma}{m}}^{-1}\,t^{-\frac{N+2\gamma}{2m}}\,\dx^{\gamma}\dy^{\gamma}
\exp\bra{-c_2 \frac{\modulus{x-y}^{\frac{2m}{2m-1}}}{t^{\frac{1}{2m-1}}}-st}$$
for some positive constants $c_1$ and $c_2$ and where $s$ is the least eigenvalue and $0\leq \gamma < m-\frac{N}{2}$
\end{theorem}
\Proof
We demonstrate the proof in two stages. Firstly optimising over $\psi$ and then optimising over
$\lambda$\\
From lemma \ref{lem:sort} we know that 
$$k\bra{t,x,y}=e^{-\lam\psi\bra{x}}k_{\lambda,\psi}\bra{t,x,y}e^{\lam\psi\bra{y}}$$
recalling the definition  of $\psi$ from (\ref{psipsi}), for some unit vector $\vecy{a}$
$$\psi\bra{\vecy{y}}-\psi\bra{\vecy{x}}=\ip{\vecy{y}-\vecy{x}}{\vecy{a}}$$
If $\vecy{a}$ is such that 
$$\psi\bra{\vecy{y}}-\psi\bra{\vecy{x}}=-|\vecy{y}-\vecy{x}|$$
then substitution into the estimate (\ref{twistedbound}) yields
\begin{equation}\label{lastlambda}
|k\bra{t,x,y}| \leq \frac{c}{\epsilon t^{1-\epsilon}}\,\dx^{\gamma}\dy^{\gamma}
e^{\small[c_2\bra{1+s}^{2m}\lambda^{2m}-s\small]t-\lam|x-y|} 
\end{equation}
Optimising the exponent of the RHS over $\lambda$ we find 
$$ \lambda^{2m-1} = \frac{\modulus{x-y}}{2mc_2\bra{1+s}^{2m}t}$$
Substituting back into (\ref{lastlambda}) we find there is a positive constant $c_2$
\begin{equation}
|k\bra{t,x,y}| 
\leq 
\frac{c}{\epsilon\ t^{1-\epsilon}}\,\dx^{\gamma}\dy^{\gamma}
\exp\bra{-c_2 \frac{\modulus{x-y}^{\frac{2m}{2m-1}}}{\bra{\bra{1+s}^{2m}t}^{\frac{1}{2m-1}}}-st} 
\end{equation}
Recalling that $\epsilon = 1-\tfrac{N+2\gamma}{2m}$  we have
the required estimate.\\
\pend


\appendix
\appendixpage

\section{Polynomial and Operator Symbol Inequalities}
When perturbing our quadratic form $Q$, the resulting twisted form $Q_{\lambda,\psi}$ generates cross terms of of the form 
$\lambda^pD^{\vecy{\alpha}}$. We can use ellipticity and Fourier transform to estimate these terms with the polyharmonic 
$\bra{-\lap}^{\frac{m}{2}}$ for $m>p+\modulus{\vecy{\alpha}}$
\begin{eqnarray*}
\norm{\lam^p D^{\vecy{\alpha}}f}^2_2 &\leq &
\hint \bra{c_1\modulus{\vecy{\xi}}^{p+\modulus{\vecy{\alpha}}}+c_2|\lam|^{p+\modulus{\vecy{\alpha}}}}|\ft|^2 d^N\vecy{\xi} \\
& \leq & c_1\bra{\mlap{\frac{m}{2}}^2_2+1}+c_2\bra{1+\lam^{2m}\norm{f}^2_2}
\end{eqnarray*}
However a more detailed decomposition of these polynomials is required to attain tighter bounds on $\norm{\lam^p D^{\vecy{\alpha}}f}^2_2$.
so that we can show exponential decay of the heat kernel in long time asymptotics.
\par

\begin{lemma}\label{lem:basic}
If $a$ , $b$ , $p$ and $q$ are all positive constants then $$a^{p}b^{q} \leq\varepsilon a^{p+q} +  c_{p,q}\varepsilon^{-\frac{p}{q}}
b^{p+q}$$ for all $\varepsilon >0$
where $c_{p,q}:=\bra{\frac{p}{p+q}}^{\frac{p}{q}}-\bra{\frac{p}{p+q}}^{1+\frac{p}{q}} $ is strictly positive and strictly less than 1.
\end{lemma}
\Proof
Follows from maximising  $x^{p}-\varepsilon{x^{p+q}}$ and substituting $x=\frac{a}{b}$. 
\pend
\begin{lemma}\label{lem:bond}
Given  $p>q>0$  there is a strictly positive constant $c_{q,p}$  such that for all
$f \in \kik{\Ome}$
$$ \mlap{q}_2 < c_{q,p}\mlap{p}_2$$
\end{lemma}
\Proof
Follows from the Spectral theorem. 
\pend
\begin{lemma}\label{lem:main}
If $f\in \kik{\Ome}$,  $\lam \in \R$ and $p$ is a positive integer
then whenever $\vecy{r}$ is a multi-index for which $|\vecy{r}|<p$ 
\begin{eqnarray*}
\norm{\lam^{p-|\vecy{r}|} D^{\vecy{r}}f}_{2}
\leq \varepsilon \mlap{\frac{p}{2}}_{2}+
\varepsilon^{-\frac{|\vecy{r}|}{p-|\vecy{r}|}}{{|\lam|}^{p}}\norm{f}_{2}
\end{eqnarray*}
for any $\varepsilon>0$
\end{lemma}
\Proof
From the isometry of the Fourier transform we see that
\begin{eqnarray*}
\norm{\lam^{p-|\vecy{r}|} D^{\vecy{r}}f}^2_2 &=& \mint \lam^{2(p-|\vecy{r}|)}\,
|D^{\vecy{r}}f |^{2} \, d^{N}\vecy{x} \\
&=& \hint \lam^{2(p-|\vecy{r}|)} { |(i\vecy{\xi})^{\vecy{r}}|}^{2}\,|\widehat{f}(\xi)|^{2} \ d^{N}\vecy{\xi}
\end{eqnarray*}
then using 
${|(i\vecy{\xi})^{\vecy{r}}|}^{2}\leq\modulus{\vecy{\xi}}^{2|\vecy{r}|}$ and lemma \ref{lem:basic}. it follows that
\begin{eqnarray*}
\norm{\lam^{p-|\vecy{r}|} D^{\vecy{r}}f}^{2}_{2} &\leq & \hint \left[\varepsilon \, \modulus{\vecy{\xi}}^{2p} +
\varepsilon^{-\frac{|\vecy{r}|}{p-|\vecy{r}|}}\, {{|\lam|}^{2p}}\,  \right]
|\widehat{f}(\xi)|^{2} \ d^{N}\underline{\xi}\\
& \leq &\varepsilon \mlap{\frac{p}{2}}_{2}^{2}+
\varepsilon^{-\frac{|\vecy{r}|}{p-|\vecy{r}|}}\, |\lam|^{2p}\,  \norm{f}_{2}^{2}		 		  		 		   
\end{eqnarray*}
Rescaling $\varepsilon$ recovers the required inequality.  \\
\pend

\begin{lemma}\label{lem:epsilon}
If $\vecy{r}$ and $\vecy{s}$ are two multi-indicies such that $|\vecy{s}| \leq p-1$ and $|\vecy{r}| \leq p$
where $p$ is a positive integer and $\lam \in \R$
then 
\begin{equation*}
\norm{\lam^{p-|\vecy{r}|}\bra{-\lap}^{\frac{|\vecy{r}|}{2}}f}_2\,
\norm{\lam^{p-|\vecy{s}|}\bra{-\lap}^{\frac{|\vecy{s}|}{2}}f}_2  
\leq \varepsilon \mlap{\frac{p}{2}}^2_2 + 2^{2p-1}\varepsilon^{1-2p}\,\lam^{2p}\,\norm{f}^2_2
\end{equation*}
for all $\varepsilon < 2$
\end{lemma}
\Proof
{\bf When $|\vecy{r}|=p$}\\
\\
We apply lemma \ref{lem:main}
\begin{eqnarray*}
\norm{\lam^{p-|\vecy{r}|}\bra{-\lap}^{\frac{|\vecy{r}|}{2}}f}_{2}
\norm{\lam^{p-|\vecy{s}|}\bra{-\lap}^{\frac{|\vecy{s}|}{2}}f}_{2}
& \leq &
\mlap{\frac{p}{2}}_{2}
\bra{\varepsilon \mlap{\frac{p}{2}}_{2} +
\varepsilon^{-\frac{|\vecy{s}|}{p-|\vecy{s}|}}
|\lam|^{p}\norm{f}_{2}} \\
& =& \epsilon\mlap{\frac{p}{2}}^{2}_{2} + \bra{\varepsilon^{-\frac{|\vecy{s}|}{p-|\vecy{s}|}}\,
|\lam|^{p}\,\norm{f}_{2}\,\mlap{\frac{p}{2}}_{2} }
\end{eqnarray*}
Using the estimate  $|ab|<\frac{1}{\varepsilon}a^{2}+\varepsilon b^{2}$
we show that the above is
\begin{eqnarray*}
& \leq & \varepsilon\mlap{\frac{p}{2}}^{2}_{2} + \frac{1}{\varepsilon}.
\varepsilon^{-\frac{2|\vecy{s}|}{p-|\vecy{s}|}}|\lam|^{2p}\norm{f}_{2}^{2}+
\varepsilon\mlap{\frac{p}{2}}^{2}_{2}  \\
& \leq & 2\varepsilon\mlap{\frac{p}{2}}^{2}_{2}+
\varepsilon^{-\frac{p+|\vecy{s}|}{p-|\vecy{s}|}}|\lam|^{2p}\norm{f}_{2}^{2}
\end{eqnarray*}
Imposing the condition $\varepsilon< 1$ and maximizing
$\varepsilon^{-\frac{p+|\vecy{s}|}{p-|\vecy{s}|}}$ over $\vecy{s}$ we have the inequality
\begin{eqnarray*}
\norm{\lam^{p-|\vecy{r}|}\bra{-\lap}^{\frac{|\vecy{r}|}{2}}f}_{2}
\norm{\lam^{p-|\vecy{s}|}\bra{-\lap}^{\frac{|\vecy{s}|}{2}}f}_{2}& \leq &
2\varepsilon\mlap{\frac{p}{2}}^{2}_{2}+
\varepsilon^{1-2p}|\lam|^{2p}\norm{f}_{2}^{2}
\end{eqnarray*}
The proof for this case is completed on rescaling $\varepsilon$.\\ \\
{\bf When $|\vecy{r}|<p$}\\
\\
Follows from lemma \ref{lem:main}.\\
\pend
\par Having completed this decomposition, we obtain an estimate for lower order operators
in terms of the higher order operator.

\begin{lemma}\label{lem:stephen}
There is a positive constant $c_1$ such for all $\rho>0$, $\theta>0$ and positive integer $p\leq m$ 
\begin{equation*}
\mlap{\frac{p}{2}}^{2}_{2}+\rho\lam^{2p} \norm{f}^{2}_{2} \leq 
c_1\bra{1+\theta}\Quad{f} + c_1\rho \bra{1+\frac{\theta s}{\rho}}^{2m} \lambda^{2m}\norm{f}^{2}_{2}
\end{equation*}

for all  $f\in \kik{\Ome}$ and with $s$ equal to the bottom eigenvalue.
\end{lemma}
\Proof

We can show by an application of lemma \ref{lem:basic} that for each $\mu >0$ this is 
\begin{equation}
\mlap{\frac{p}{2}}^{2}_{2}+\rho\lam^{2p} \norm{f}^{2}_{2}\leq \mlap{\frac{p}{2}}^{2}_{2}+\rho \mu \norm{f}^{2}_{2}+ \rho\mu^{-\frac{m-p}{p}} \lam^{2m} \norm{f}^{2}_{2}
\end{equation}{lem:sort}
Applying lemma \ref{lem:bond} and ellipticity whilst recalling that we have  positive-definite spectral gap $s$
\begin{eqnarray*}
\mlap{\frac{p}{2}}^{2}_{2}+\rho\lam^{2p} \norm{f}^{2}_{2}&\leq& c_1\mlap{\frac{m}{2}}^{2}_{2}+ \rho m\mu\norm{f}^{2}_{2} + \rho m\bra{1+\mu}^{m-1} \lam^{2m} \norm{f}^{2}_{2}\\
&\leq &c_1\bra{1+\theta}\Quad{f}- c_1\theta s\norm{f}^{2}_{2}+ \rho m\mu\norm{f}^{2}_{2} + \rho m\bra{1+\mu}^{m-1} \lam^{2m} \norm{f}^{2}_{2}\\
\end{eqnarray*}
we can then set $\mu=\frac{c_1\theta s}{m\rho}$\\
\pend

\section{Davies' Twisted Operators and the Canonical Functional Calculus}

The two definitions are consistent
\begin{eqnarray*}
Q_{\lam \psi}\bra{f}&=& \langle \HH^\frac{1}{2}\twist f,\HH^\frac{1}{2}\itwist f \rangle\\
&=& \langle \itwist \HH \twist f ,f\rangle\\
&=& \langle  \HH_{\lambda\psi} f ,f\rangle\\
\end{eqnarray*}
\begin{lemma}
$\HH$ and $\HH_{\lambda\psi}$ have the same spectrum.
\end{lemma}
\Proof
Given a sequence of functions $f_n$ where $\bra{\HH-z} f_n\rightarrow 0$ we have
$\bra{\HH_{\lambda\psi}-z}\bra{\itwist f_n}\rightarrow 0$\\
\pend

\begin{lemma}
For each $f\in L^2\bra{\Ome}$ and $z$ in the resolvent set of $\HH$ we have
$$\bra{z-\HH_{\lambda\psi}}^{-1}f=\itwist \bra{z-\HH}^{-1} \twist f$$
\end{lemma}
\Proof
Let $\BB=\bra{z-\HH_{\lambda\psi}}^{-1}$ then
$$\BB\bra{z-\HH_{\lambda\psi}}f=f$$
which is just
$$\BB\itwist\bra{z-\HH}\twist f=f$$
and out statement follows.\\
\pend
\begin{corollary}
For f $\in C_0\bra{R}$ we have a canonical functional calculus for the twisted operator given
by $$ f\bra{\HH_{\lambda\psi}} :\ =\ \itwist f\bra{\HH} \twist$$
\end{corollary}

\begin{lemma}\label{lem:sort}
If $f \in C_0\bra{R}$ and $k_f\bra{x,y}$ is the integral kernel of the operator
$f\bra{\HH}$ then $f\bra{\HH_{\lambda\psi}}$ has integral kernel $k_{f,{\lambda\psi}}\bra{x,y}$
where
$$k_{f,{\lambda\psi}}\bra{x,y}=e^{-\lam\psi\bra{x}}k_f\bra{x,y}e^{\lam\psi\bra{y}}$$
\end{lemma}
\Proof
Let $v \in L^2\bra{\Ome}$ then 
\begin{eqnarray*}
\mint k_{f,{\lambda\psi}}\bra{x,y} \ v\bra{y} dy &=& \bra{f\bra{\HH_{\lambda\psi}}v}\bra{x}\\
&=&\bra{\itwist f\bra{\HH}\twist v}\bra{x}\\
&=& e^{-\lambda \psi \bra{x}}\mint k_f\bra{x,y}e^{\lambda \psi \bra{y}}v\bra{y} dy \\
\end{eqnarray*}
\pend

 \begin{center}
 {\Large{\bf{ Acknowledgements}}}
 \end{center}
This research was funded by an EPSRC Ph.D grant 95-98 at Kings College, London. I would like to thank E.Brian Davies
for giving me this problem and his encouragement since. I am very grateful to Gerassimos Barbatis for all
his helpful comments and his prodding me to complete this paper. I am indebted to Anita for all her support.
 \par
Narinder Claire \newline
Global Equities \& Commodity Derivatives Quantitative Research \newline
BNP Paribas London \newline
10 Harewood Avenue\newline
London\newline
NW1 6AA \newline
e-mail: narinder.claire@uk.bnpparibas.com
\end{document}